\documentstyle[a4,12pt]{article}
\begin{document}
\begin{center}
{\bf INTEGRABILITY OF THE GAUSS-CODAZZI-MAINARDI EQUATION IN 2+1
DIMENSIONS}
\end{center}
\begin{center}
{\bf R.Myrzakulov,  Kur. Myrzakul}
\end{center}
\begin{center}
 Institute of Physics and Technology, 480082  Alma-Ata-82, Kazakhstan
\end{center}
\section{Introduction}
At present, it is clear that a lot of the results of modern 
soliton theory
can be found in the old text-books on differential geometry and
numerous old papers of XIX century (see,
for instance, [9] and references therein).    
So that one of important root of
soliton theory lies in differential geometry (see, e.g. 
[1-7] and references therein).  
Itself soliton equations appeared as the
some particular cases of the Gauss-Codazzi-Mainardi 
equations (GCME). At the same time, so-called integrable 
spin systems are some reductions of the Gauss-Weingarten 
equation (GWE) [3]. In
modern soliton therminology this means that the GCME and 
the GWE admit some integrable cases. On the other hand, as shown 
quite  resently
by V.E.Zakharov in the remarkable paper [5], 
the (1+1)-dimensional  GCME 
is integrable. The aim of this
notice is to consider the  integrability aspects of the GCME 
in 2+1 dimensions and its connections with the other integrable 
equations such as the Yang-Mills-Higgs-Bogomolny equation (YMHBE)
and the self-dual Yang-Mills equation (SDYME).

\section{The GWE  and GCME in 1+1 dimensions}

Let us consider the GWE
$$
r_{xx}=\Gamma^1_{11} r_x+\Gamma^2_{11} r_t +Ln \eqno(1a) $$
$$ r_{xt}=\Gamma^1_{12}+\Gamma^2_{12} r_t+M n \eqno(1b) $$
$$ r_{tt}=\Gamma^1_{22} r_x+\Gamma^2_{22} r_t +Nn \eqno(1c)$$
$$ n_x=P_{11} r_x+ P_{12} r_t \eqno(1d)$$
$$ n_t=P_{21} r_x+P_{22} r_t. \eqno(1e) $$
If we introduce the vector
 $Z=(r_x,r_t, n)^t, $ then the GWE takes the form
$$ Z_x=A^{\prime}Z, \quad Z_t=C^{\prime}Z \eqno(2) $$
where
$$
A^{\prime}=\left( \begin{array}{ccc}
\Gamma^1_{11} & \Gamma^2_{11} & L \\
\Gamma^1_{12} & \Gamma^2_{21} & M \\
P_{11} & P_{12} & 0
\end{array} \right), \quad
C^{\prime}=\left (\begin{array}{ccc}
\Gamma^1_{12} & \Gamma^2_{12} & M \\
\Gamma^1_{22} & \Gamma^1_{22} & N \\
P_{12} & P_{22} & 0
\end{array} \right). \eqno(3)
$$
Let us introduce the orthogonal basis
$$
{\bf  e}_1=\frac{{\bf r}_{x}}{\sqrt{E}},\quad {\bf e}_2=
{\bf n}, \quad {\bf e}_3={\bf e}_1 \wedge {\bf e}_2 \eqno(4) $$
and let  ${\bf  e}^{2}_1=\beta=\pm 1, \quad  {\bf e}^{2}_2=
{\bf e}^{2}_3=1$.
Then the GWE (1) takes the form
$$
\left (\begin{array}{c}
e_1 \\ e_2 \\ e_3
\end{array} \right)_x= A
\left (\begin{array}{c}
e_1 \\ e_2 \\ e_3
\end{array} \right), \quad
\left (\begin{array}{c}
e_1 \\ e_2 \\ e_3
\end{array} \right)_t= C
\left(\begin{array}{c}
e_1 \\ e_2 \\ e_3
\end{array} \right) \eqno(5)
$$
where (about our notation see, e.g. [3])
$$  A=\left(\begin{array}{ccc}
0 & k & -\sigma \\
-\beta k & 0 & \tau \\
\beta\sigma & -\tau & 0
\end{array} \right), \quad
C=
\left(\begin{array}{ccc}
0 & \omega_{3} & -\omega_{2} \\
-\beta\omega_{3} & 0 & \omega_{1} \\
\beta\omega_{2} & -\omega_{1} & 0
\end{array} \right).
 \eqno(6)
$$
Here $\beta ={\bf e}^{2}_{1}=\pm 1, \quad {\bf e}^{2}_{2}={\bf e}^{2}_{3}=
1$.
Below we consider the case $\beta=+1$.
Integrability conditions for the GWE (5) give the  following
GCME
$$ A_t-C_x+[A,C]=0 \eqno(7) $$
or in elements
$$
k_t-\omega_{3x}=\tau\omega_2-\sigma\omega_1 \eqno (8a)
$$
$$
\tau_t-\omega_{1x}=\sigma\omega_3- k\omega_2 \eqno (8b)
$$
$$
\sigma_t-\omega_{2x}=k\omega_1-\tau\omega_3. \eqno (8c)
$$
Note that the GWE (5) can be rewritten in the 2$\times$2 matrix
form as
$$ g_x=Ug,  \quad  g_t=Wg \eqno(9) $$
with
$$ U=
\frac{1}{2i}\left ( \begin{array}{cc}
\tau & k+i\sigma \\
k-i\sigma & -\tau
\end{array} \right),
W=
\frac{1}{2i}\left ( \begin{array}{cc}
\omega_1 & \omega_3+i\omega_2 \\
\omega_3-i\omega_2 & -\omega_1
\end{array} \right).   \eqno(10)
$$
Hence we get the following form of the  GCME
$$U_t-W_x+[U,W]=0. \eqno(11) $$
We note that the system (9) can be considered as the Lax 
representation (LR) for the GCME (11).

\section{The GWE and GCME  in 2+1 dimensions}

Now we consider the 3-dimensional manifold $V^{3}$ immersed 
into the $R^{n}$. The corresponding (2+1)-dimensional GWE is 
given by
$$ \left(\begin{array}{c}
e_1 \\ e_2 \\ \vdots \\ e_n
\end{array} \right)_x=A
\left(\begin{array}{c}
e_1 \\ e_2 \\ \vdots \\e_n
\end{array} \right),
\left(\begin{array}{c}
e_1 \\e_2 \\ \vdots \\ e_n
\end{array} \right)_y=B
\left(\begin{array}{c}
e_1 \\e_2 \\ \vdots \\e_n
\end{array} \right),
\left(\begin{array}{c}
e_1 \\e_2 \\  \vdots \\ e_n
\end{array} \right)_t=C
\left(\begin{array}{c}
e_1 \\e_2 \\  \vdots \\ e_n
\end{array} \right)
\eqno(12)
$$
where   $A, B, C\in so(l,m), l+m=n$.
The compatibility conditions of these equations give the 
 following system
$$
A_y-B_x +[A, B]=0 \eqno(13a)
$$
$$
 A_t-C_x+[A, C]=0 \eqno(13b)
$$
$$
B_t-C_y+[B, C]=0. \eqno(13c)
$$
Usually this system is called the (2+1)-dimensional
GCME. In (13), e.g. 
for the case $n=3$, the matrix $B$ has the form
$$ B=
\left(\begin{array}{ccc}
0 & m_{3} & -m_{2} \\
-\beta m_{3} & 0 & m_{1} \\
\beta m_{2} & -m_{1} & 0
\end{array} \right).
 \eqno(14)
$$

\subsection{Lax representation for the (2+1)-dimensional GCME}
The GCME (13) admits two types LR.
\subsubsection{LR of 1-type}
The LR of the GCME (13) has the form
$$ g_x=Ug, \quad g_y=Vg, \quad  g_t=Wg \eqno(15) $$
with (for the case $n=3$) $ U, W$ given by (10) and
$$
V=
\frac{1}{2i}\left ( \begin{array}{cc}
m_1 & m_3+im_2 \\
m_3-im_2 & -m_1
\end{array} \right). \eqno(16)
$$
Hence we get the new form of the (2+1)-dimensional GCME
$$U_t-W_x+[U,W]=0 \eqno(17a) $$
$$ U_y-V_x+[U,V]=0 \eqno(17b) $$
$$ V_t-W_y+[V,W]=0. \eqno(17c) $$
Using the transformation
$$
\psi=ge^{I_{1}x+I_{2}y+I_{3}t} \eqno(18)
$$
the LR (15) can be rewritten in the following form
$$
\psi_x= A\psi-\psi I_{1}(\lambda), \quad 
\psi_y= B\psi-\psi I_{2}(\lambda), \quad
\psi_t= C\psi-\psi I_{3}(\lambda). \eqno(19)
$$
Here $I_{k}$ are the some  constant diagonal matrices. 
\subsubsection{LR of  2-type}
There exist the other form of the LR for the GCME (13)
which reads as
$$
L_{1}\Psi=(A_{t}+ A_{y}-\lambda A_{x})\Psi, \quad
L_{2}\Psi=(-\lambda A_{t}+ \lambda A_{y}+A_{x})\Psi
\eqno(20)
$$
where 
$$
L_{1}=-\partial_{t}-\partial_{y}+\lambda \partial_{x},  \quad
L_{2}=\lambda\partial_{t}-\lambda\partial_{y}- \partial_{x}. 
\eqno(21)
$$ 
So the GCME (13) can be integrated by using the LR (19) or (20).
\section{The (2+1)-dimensional GCME as the particular
case of the Yang-Mills-Higgs-Bogomolny  equation}
It is interesting to note that the (2+1)-dimensional GCMRE (13) 
is the
particular
case of the  YMHBE. In the simplest  case, the YMHBE reads as
$$
\Psi_{t}+[\Psi, C]+A_y-B_x +[A, B]=0, \eqno(22a)
$$
$$
\Psi_{y}+[\Psi, B]+ A_t-C_x+[A, C]=0, \eqno(22b)
$$
$$
\Psi_{x}+[\Psi, A]+B_t-C_y+[B, C]=0. \eqno(22c)
$$
Hence  as $\Psi=0$ we get the GCME (13). 
The YMHBE is integrable. The corresponding
LR, e.g,  has the form
$$
L_{1}\Psi=(C+ B-\lambda A-\lambda \Phi)\Psi, \quad
L_{2}\Psi=(-\lambda C+ \lambda B+A- \Phi)\Psi
\eqno(23)
$$
where $L_{k}$ are given by (21).
The covariant form of the LR (23) is
$$
(D_{t}+D_{y}-\lambda D_{x})\Psi=\lambda \Phi\Psi,  \quad
\lambda D_{t}-\Lambda D_{y}-D_{x})\Psi=- \Phi\Psi
\eqno(24)
$$
where 
$$
D_{i}=\partial_{i}\Phi+[A_{i}, \Phi] \eqno(25)
$$
and $A_{t}=A, A_{y}=B, A_{x}=A$.
\section{The (2+1)-dimensional GCME as the particular
case of the  self-dual Yang-Mills equation}

Finally we observe that 
the GCME (13) is also the particular case of the SDYME. 
The SDYME reads as 
$$
F_{\alpha\beta}=0,  \quad F_{\bar\alpha\bar\beta}=0,  \quad
F_{\alpha\bar\alpha}+
F_{\beta\bar\beta}=0 \eqno(26)
$$
where
$$
F_{\mu\nu}=A_{\nu,\mu}-A_{\mu,\nu}-[A_{\mu},A_{\nu}].  \eqno(27)
$$
As known the LR for the SDYME is given by [8, 10]
$$
(\partial_{\alpha}+\lambda\partial_{\bar\beta})\Psi=
(A_{\alpha}+\lambda A_{\bar\beta})\Psi, \quad
(\partial_{\beta}-\lambda\partial_{\bar\alpha})\Psi=
(A_{\alpha}+\lambda A_{\bar\beta})\Psi
\eqno(28)
$$
where $\lambda$ is the spectral parameter.
If in the SDYME (26) we take
$$
A_{\alpha}=-iC, \quad A_{\beta}=A-iB, 
\quad A_{\bar\beta}=A+iB, 
\quad A_{\bar\alpha}=iC \eqno(29)
$$ 
and assume that $A,B,C$ are independent of $z$ then we 
obtain the GCME (13).

\section{Conclusion}
In this notice the integrability of the (2+1)-dimensional 
GCME is discussed. It is shown that this equation is 
integrable in the
sense that for it the corresponding LR with the spectral parameter 
exist. Also we have established that the GCME in 2+1 dimensions 
is the 
particular case
of the YMHBE and of the SDYME.
\section{Akcnowledgment}
R.M. thanks the Prof. Antoni Sym for sending of the book [7].
This work was supported by INTAS, grant 99-1782. 

\end{document}